\documentclass[11pt]{amsart}
\usepackage[a4paper, margin=0.5in]{geometry}
\usepackage[all]{xy}
\usepackage{pdfsync}
 \usepackage[colorlinks=false]{hyperref}
\usepackage{amssymb}
\usepackage{amsmath}
\usepackage{autonum}
\usepackage{graphicx}
\usepackage[english]{babel}
\usepackage{epsfig}
\usepackage{multirow}
\usepackage{color}
\usepackage{enumerate}
\usepackage[T1]{fontenc}
\usepackage{aurical, pbsi,calligra}
\swapnumbers
\usepackage{amsmath,systeme}
\usepackage{amsmath} 
\usepackage{graphics}
\usepackage{rotating}
\usepackage{lipsum}
\usepackage{lscape}
\usepackage{tabularx}
\usepackage{afterpage}
\usepackage{capt-of}% or use the larger `caption` package

\usepackage{booktabs}
\usepackage{siunitx}

\usepackage{lipsum} % For dummy text
\usepackage{pdflscape} % For landscape tables
\usepackage{geometry} % For page dimensions
\usepackage{tabularx} % For table width control

\theoremstyle{definition}
\newtheorem{definition}{Definition}[section]

\theoremstyle{plain}

\newtheorem{theorem}[definition]{Theorem}
\usepackage[all]{xy}
\usepackage{pdfsync}
\usepackage{ytableau}

\newrobustcmd*{\mytriangle}[1]{\tikz{\filldraw[draw=#1,fill=#1] (0,0) --
(0.2cm,0) -- (0.1cm,0.2cm);}}

\begin{document}
\title [A Collocation Heat Polynomials Method for Inverse Stefan Problems] {A Collocation Heat Polynomials Method for One-Dimensional Inverse Stefan Problems}

\author[O. Narbek]{Orazbek Narbek}

\address{SDU University, Department of Mathematics, Almaty, Kazakhstan}

\email{narbekov.o@gmail.com}

\author[S. A. Kassabek]{Samat A. Kassabek}

\address{Astana IT university, Department of Computational and Data Science, Astana, Kazakhstan}

\email{samat.kassabek@astanait.edu.kz}

\author[T. Nauryz]{Targyn Nauryz}

\address{Kazakh British Technical University, Institute of Mathematics and Mathematical Modeling, Narxoz University, Almaty, Kazakhstan}

\email{targyn.nauryz@gmail.com }

\subjclass[2010] {80A22, 80A23, 80M30, 35C11}

\keywords{Inverse Stefan problems, approximate solution, heat polynomials method, moving boundary.}

\thanks{The authors were supported by the grant AP19675480 "Problems of heat conduction with a free boundary arising in modeling of switching processes in electrical devices".}

\maketitle

%\sloppy

\bigskip

\begin{abstract} 
The inverse one-phase Stefan problem in one dimension, aimed at identifying the unknown time-dependent heat flux \( P(t) \) with a known moving boundary position \( s(t) \), is investigated. A previous study \cite{Samat} attempted to reconstruct the unknown heat flux \( P(t) \) using the Variational Heat Polynomials Method (VHPM). In this paper, we develop the Collocation Heat Polynomials Method (CHPM) for the reconstruction of the time-dependent heat flux \( P(t) \). This method constructs an approximate solution as a linear combination of heat polynomials, which satisfies the heat equation, with the coefficients determined using the collocation method. To address the resulting ill-posed problem, Tikhonov regularization is applied. As an application, we demonstrate the effectiveness of the method on benchmark problems. Numerical results show that the proposed method accurately reconstructs the time-dependent heat flux \( P(t) \), even in the presence of significant noise. The results are also compared with those obtained in \cite{Samat} using the VHPM.
\end{abstract}

\bigskip

\section{\textbf{Introduction}}
The heat polynomials \(v_n(x, t)\) were introduced by Appell in 1892 \cite{Appell} as the coefficients of the generating function:
$$
g(x, t, z) = e^{x z + t a^2 z^2},
$$
where the function \(g(x, t, z)\) satisfies the heat equation \(g_t - a^2 g_{xx} = 0\) for any value of \(z\). This function can be expanded in powers of \(z\) using a Taylor series:
$$
g(x, t, z) = \sum_{n=0}^{\infty} v_n(x, t) z^n.
$$
Rosenbloom and Widder, in their series of works \cite{Rosenbloom}, \cite{Widder} - \cite{ Widder4}, along with Haimo \cite{Haimo}, investigated the expansion of solutions \(u(x, t)\) to the heat equation
\begin{equation}
    \frac{\partial u}{\partial t} = a^2 \frac{\partial^2 u}{\partial x^2}, \label{14}
\end{equation}
using heat polynomials \(v_n(x, t)\). These polynomials are given by the expression:
\begin{equation}
    v_n(x, t) = \sum_{m=0}^{\left\lfloor n / 2 \right\rfloor} \frac{a^{2m} n!}{m! (n - 2m)!} x^{n - 2m} t^m. \label{08}
\end{equation}
It is evident that this linear combination of heat polynomials also satisfies the heat equation (\ref{14}). 

The core idea behind this method is to construct exact analytical or approximate solutions to the problem as a linear combination of the polynomials in (\ref{08}), much like in the Trefftz method \cite{Trefftz}, where unknown coefficients are determined from initial and boundary conditions. This approach has been further developed in the works of Yano \cite{Yano}, Futakiewicz \cite{Futakiewicz}, Hożejowski \cite{Hozejowski}, and Grysa \cite{Grysa}.

Numerous studies have explored the application of the heat polynomial method (HPM) and the Trefftz method for various heat conduction problems on fixed domains, see e.g., \cite{Grysa1} - \cite{Grysa3}, \cite{Liu} - \cite{Maciejewska2}. However, the literature on applying HPM to heat problems with moving boundaries is relatively scarce; see \cite{Samat} - \cite{Samat4}, \cite{Kharin3, Kharin4, Reemtsen}.

Heat conduction problems with moving boundaries and phase transitions, known as Stefan problems, arise in various real-world and engineering applications, such as melting, welding, freezing, and cooling \cite{Gupta, Rubinstain}. The direct Stefan problem involves determining the temperature distribution and the position of a moving boundary over time, given the heat flux or boundary conditions. This problem typically models processes where phase change, such as melting or evaporation, drives the movement of the boundary. In contrast, inverse Stefan problem focuses on identifying unknown parameters, such as the time-dependent heat flux or the boundary's position, based on known temperature data or other boundary conditions \cite{Goldman}. Due to its ill-posed nature, this problem is more complex and commonly encountered in situations where direct measurements are challenging. An example of a real-world application of the inverse Stefan problem is in determining the heat flux components in low-voltage electric arc systems \cite{Kharin2}. This is particularly relevant since experimental methods for measuring heat fluxes can often result in significant errors, making the use of inverse methods highly valuable.

Many studies have addressed inverse Stefan problems \cite{Chifaa}, \cite{Grzymkowski1}, \cite{Johansson1} - \cite{Slota}, \cite{Rad}, \cite{Reemtsen}. One effective method for solving these problems is the heat potential method, which reduces boundary value problems to integral equations. However, when the domain is degenerate at the initial time, singularities in these integral equations introduce additional difficulties. Numerical methods like finite difference and finite element techniques \cite{Smith, Zabaras} are widely used but can be computationally expensive, especially when fine discretization is required. To overcome these limitations, meshless methods \cite{Johansson} and method of fundamental solutions \cite{Johansson1} have been developed. While flexible, these methods often face stability and convergence issues, particularly with noisy data. The slow convergence of many iterative methods and the need for accurate solutions with low computational cost have motivated this study. We aim to improve the efficiency and stability of inverse one-phase Stefan problem solutions using the CHPM, which addresses several of these limitations.

The collocation method has gained prominence as an effective and versatile tool for solving inverse heat conduction problems (IHCPs) and Stefan problems. Unlike traditional mesh-based methods, the collocation method employs radial basis functions or other approximating functions to solve partial differential equations without the necessity of a predefined mesh. This meshless approach provides significant flexibility, especially when dealing with complex geometries and irregular domains, while simultaneously reducing computational costs associated with mesh generation \cite{Khan}. Over time, the collocation method has evolved to include various regularization techniques aimed at addressing the ill-conditioning that typically arises in inverse problems, particularly IHCPs, where small errors in data can lead to large errors in solutions \cite{Shahnazari}. For example, recent advancements such as the local meshless methods \cite{Wang} and physics-informed neural networks (PINNs) \cite{Zhang} have further extended the capabilities of the collocation method for complex multi-layered media and nonlinear heat conduction problems, offering promising solutions to both forward and inverse problems. These developments have significantly enhanced the stability and accuracy of the method, making it an increasingly popular approach for tackling one-phase Stefan problems \cite{Khan, Shahnazari, Vrankar}.

In this article, we present the approximate solutions of one-phase inverse Stefan problems, focusing on the development of a novel collocation heat polynomials method. The primary contribution of this work lies in the application of CHPM for solving inverse Stefan problems, offering an efficient and accurate alternative to existing methods. Our goal in this series of studies \cite{Samat} - \cite{Samat4} is to create and refine an approximate method that can effectively handle Stefan-type problems. To demonstrate the method's efficacy and investigate the convergence of the approximate solutions, we apply CHPM to two numerical examples, illustrating its potential as a robust tool for inverse Stefan problems.

The paper is organized as follows. Formulation of the one-dimensional inverse Stefan problems are presented in Section 2. In Section 3, we give a brief sketch of the heat polynomials method to solve the Stefan type problems. In Section 4, we consider an application of the CHPM to two numerical examples. In order to investigate the stability of the approximate solution, random additive noise is added to the input data. Summary and Concluding remarks are discussed briefly in Section 5.

\section{\textbf{One-phase Stefan Problem}}

Consider the following one-dimensional, one-phase Stefan problem:
\begin{equation}
    \frac{\partial u}{\partial t} = a^2 \frac{\partial^2 u}{\partial x^2} \quad \text{in} \quad 0 < x < s(t), \quad 0 < t \leq T, \label{02}
\end{equation}
where \(a\) represents the thermal diffusivity, and \(s(t) > 0\) is a free boundary that is assumed to be a sufficiently smooth function. The heat conduction domain is defined as \(D = (0, s(t)) \times (0, T)\), with its closure given by \(\overline{D} = [0, s(t)] \times [0, T]\). 

The initial condition is provided as:
\begin{equation}
    u(x, 0) = f(x), \quad 0 \leq x \leq s(0), \label{03}
\end{equation}
with Neumann boundary conditions specified as:
\begin{equation}
    - \lambda \frac{\partial u(0, t)}{\partial x} = P(t), \quad 0 < t \leq T, \label{04}
\end{equation}
and Dirichlet and Stefan conditions imposed on the moving boundary \(x = s(t)\):
\begin{equation}
    u(s(t), t) = u^*, \quad 0 < t \leq T, \label{05}
\end{equation}
\begin{equation}
    - \lambda \frac{\partial u(x, t)}{\partial x} \Bigr|_{x = s(t)} = L \gamma \frac{ds(t)}{dt}, \quad 0 < t \leq T, \label{96}
\end{equation}
along with the initial condition on the moving boundary:
\begin{equation}
    s(0) = s_0. \label{93}
\end{equation}

For the direct one-phase Stefan problem, we seek to determine the free boundary \(x = s(t) \in C^1([0, T])\) and the unknown temperature distribution \(u(x, t) \in C^2(D) \cap C^1(\overline{D})\), which satisfy Eqs. \((\ref{02})\), \((\ref{03})\), \((\ref{04})\), \((\ref{05})\), \((\ref{96})\), and \((\ref{93})\). The parameters \(a\), \(\lambda\), \(u^*\), \(L\), \(\gamma\), \(f(x)\), \(P(t)\), and \(s_0\) are assumed to be known.

Unlike the direct Stefan problem, numerous variations of inverse Stefan problems can be defined. In this paper, we focus on investigating a classical inverse Stefan problem where the Neumann boundary condition at \(x=0\) is not prescribed and must be reconstructed, along with the temperature distribution, assuming the position of the moving boundary \(s(t)\) is known a priori.

\section{\textbf{The Collocation Heat Polynomials Method}}

In this study, we employ the collocation heat polynomials method to solve the inverse one-phase Stefan problem. This method approximates the temperature distribution and heat flux at the moving boundary using a linear combination of heat polynomials, which inherently satisfy the heat equation. The Stefan problem is transformed into a system of algebraic equations by applying the boundary and initial conditions through integration.

We approximate the temperature distribution \( u(x, t) \) as a linear combination of heat polynomials:
\[
    u(x, t) = \sum_{n=0}^{N} c_n v_n(x, t), \quad 0 \leq x \leq s(t), \quad 0 \leq t \leq T \label{74}
\]
where \( v_n(x, t) \) represents the \(n^{th}\)-order heat polynomials. These polynomials satisfy the heat equation, and the coefficients \(c_n\) are unknown and determined through the collocation process.

To determine these coefficients, we apply the collocation method by enforcing the boundary and initial conditions through an integration approach over finite intervals. Instead of evaluating the conditions at discrete points, we integrate the residuals over the intervals, making the system more stable for ill-posed problems.

At the moving boundary \(x = s(t)\), we integrate the Dirichlet and Neumann conditions over time intervals \(\left( t_{i-1}, t_i \right]\), for \( i = 1, 2, \dots, M \), as follows:
\[ \int_{t_{i-1}}^{t_i} \left[ u(s(t), t) - u^* \right] \, dt = \int_{t_{i-1}}^{t_i} \left[ \sum_{n=0}^{N} c_n v_n(s(t), t) - u^* \right] \, dt = 0,\] and for the Neumann condition, we have:
\[
    \int_{t_{i-1}}^{t_i} \left[ -\lambda \frac{\partial u}{\partial x} \Bigg|_{x = s(t)} - L \gamma \frac{ds(t)}{dt} \right] dt = \int_{t_{i-1}}^{t_i} \left[ -\sum_{n=0}^{N} c_n \lambda \frac{\partial v_n(x, t)}{\partial x} \Bigg|_{x = s(t)} - L \gamma \frac{ds(t)}{dt} \right] dt = 0.
\]
The initial condition \( u(x, 0) = f(x) \) is handled similarly. We divide the spatial domain \([0, s(0)]\) into \( N_x \) intervals and integrate the residuals over each interval:
\[
    \int_{x_{j-1}}^{x_j} \left[ u(x, 0) - f(x) \right] dx = \int_{x_{j-1}}^{x_j} \left[ \sum_{n=0}^{N} c_n v_n(x, 0) - f(x) \right] dx = 0.
\]
These integrals result in a system of \(N \times N\) algebraic equations that can be written in matrix form as:
\[
    \mathbf{A} \mathbf{c} = \mathbf{b},  \label{10}
\]
where \( \mathbf{A} \) contains coefficients derived from the integrals, \( \mathbf{c} \) is the vector of unknown coefficients \(c_n\), and \( \mathbf{b} \) contains the results from the integrated boundary and initial conditions.

The system of equations obtained is often ill-conditioned, meaning that small perturbations in the data can lead to significant errors in the solution \cite{Ang}. To address this, we apply Tikhonov regularization, modifying the system of equations to:
\[
    (\mathbf{A}^{\text{T}} \mathbf{A} + \beta I) \mathbf{c} = \mathbf{A}^{\text{T}} \mathbf{b},
\]
where \( \mathbf{A}^{\text{T}} \) is the transpose of \( \mathbf{A} \), \( I \) is the identity matrix, and \( \beta \) is the regularization parameter. This regularization helps to stabilize the solution by controlling the balance between fitting the data and maintaining numerical stability. In our numerical examples, we explore the effect of varying the regularization parameter \( \beta \) on the accuracy of the reconstructed heat flux.

\begin{figure}[h!]
\centering
\includegraphics[width = 16cm]{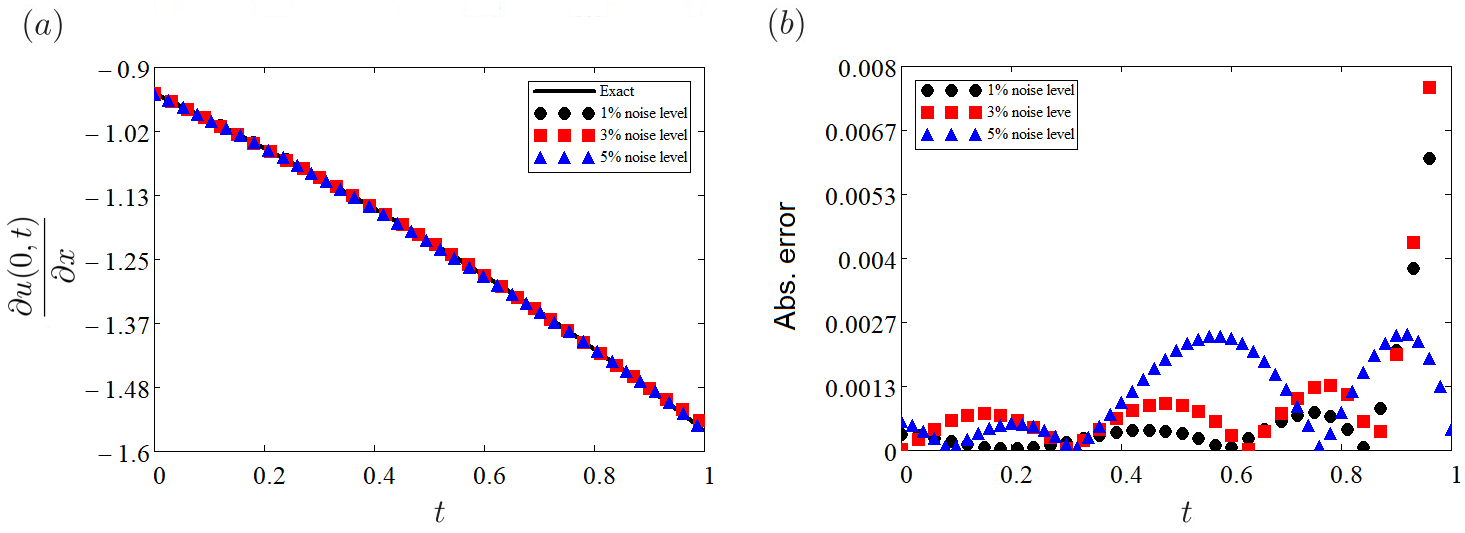}
\caption{(a) For Example 4.1 of the inverse Stefan problem (a) depicts the comparison between the exact and approximate solutions of $u_x(0, t)$, (b) illustrates the absolute error between the exact solution and $u_x(0, t)$ for different noise levels for $N=12$.}
\label{fig6}
\end{figure}

\begin{table}[ht]
\caption{The values $\Delta P$ of the relative error depending on the parameters $N$ and $\beta$ with zero noise, Example 4.1.} \centering
\scriptsize
\label{table 1}
\renewcommand{\arraystretch}{1.2} % Increases the height of each row
\resizebox{\columnwidth}{!}{
\begin{tabular}{|c|c|c|c|c|c|c|c|c|c|}
\hline 
$\beta \backslash N$  & 4 & 6 & 8 & 10 & 12 & 14 & 16 & 18  & 20 \\
\hline
0 & 0.018 & $1.90 \cdot 10^{-3}$ & $1.08\cdot 10^{-4}$ & $5.06\cdot 10^{-6}$ & $1.71 \cdot 10^{-7}$ & $5.06 \cdot 10^{-9}$ & $6.57 \cdot 10^{-11}$ & $1.67 \cdot 10^{-11}$ & $2.05 \cdot 10^{-11}$ \\
\hline 
$1 \cdot 10^{-13}$  & 0.018 & $1.91 \cdot 10^{-3}$ & $1.08 \cdot 10^{-4}$ & $5.36 \cdot 10^{-6}$ & $2.76 \cdot 10^{-6}$ & $2.61 \cdot 10^{-6}$ & $3.08 \cdot 10^{-6}$ & $2.40 \cdot 10^{-6}$ & $2.97 \cdot 10^{-6}$  \\
\hline 
$1 \cdot 10^{-12}$  & 0.018 & $1.91 \cdot 10^{-3}$ & $1.08 \cdot 10^{-4}$ & $9.8 \cdot 10^{-6}$ & $7.33 \cdot 10^{-6}$ & $7.08 \cdot 10^{-6}$ & $7.44 \cdot 10^{-6}$ & $6.93 \cdot 10^{-6}$ & $7.16 \cdot 10^{-6}$\\
\hline 
$1 \cdot 10^{-11}$  & 0.018 & $1.91 \cdot 10^{-3}$ & $1.10 \cdot 10^{-4}$ & $2.85 \cdot 10^{-5}$ & $2.54 \cdot 10^{-5}$ & $2.16 \cdot 10^{-5}$ & $2.03 \cdot 10^{-5}$ & $2.02 \cdot 10^{-5}$ & $2.06 \cdot 10^{-5}$ \\
\hline 
$1 \cdot 10^{-10}$ & 0.018 & $1.91 \cdot 10^{-3}$ & $1.23 \cdot 10^{-4}$ & $6.69 \cdot 10^{-5}$ & $6.42 \cdot 10^{-5}$ & $6.28 \cdot 10^{-5}$ & $7.28 \cdot 10^{-5}$ & $6.28 \cdot 10^{-5}$ & $7.40 \cdot 10^{-5}$  \\
\hline 
$1 \cdot 10^{-9}$  & 0.018 & $1.92 \cdot 10^{-3}$ & $2.28 \cdot 10^{-4}$ & $1.85\cdot 10^{-4}$ & $1.76 \cdot 10^{-4}$ & $1.75 \cdot 10^{-4}$ & $1.86 \cdot 10^{-4}$ & $1.79 \cdot 10^{-4}$ & $1.89 \cdot 10^{-4}$ \\
\hline 
$1 \cdot 10^{-8}$  & 0.018 & $2.01 \cdot 10^{-3}$ & $5.84 \cdot 10^{-4}$ & $5.75 \cdot 10^{-4}$ & $5.16 \cdot 10^{-4}$ & $5.14 \cdot 10^{-4}$ & $5.29 \cdot 10^{-4}$ & $5.32 \cdot 10^{-4}$ & $5.59 \cdot 10^{-4}$ \\
\hline 
$1 \cdot 10^{-7}$ & 0.018 & $2.74 \cdot 10^{-3}$ & $1.57 \cdot 10^{-3}$ & $1.53 \cdot 10^{-3}$ & $1.57 \cdot 10^{-3}$ & $1.62 \cdot 10^{-3}$ & $1.79 \cdot 10^{-3}$ & $1.69 \cdot 10^{-3}$ & $1.87 \cdot 10^{-3}$ \\
\hline 
$1 \cdot 10^{-6}$ & 0.016 & $5.45 \cdot 10^{-3}$ & $4.61 \cdot 10^{-3}$ & $4 \cdot 10^{-3}$ & $4.1 \cdot 10^{-3}$ & $4.25 \cdot 10^{-3}$ & $4.49 \cdot 10^{-3}$ & $4.5 \cdot 10^{-3}$ & $4.75 \cdot 10^{-3}$ \\
\hline 
$1 \cdot 10^{-5}$ & 0.009 & 0.011 & 0.011 & 0.011 & 0.012 & 0.013 & 0.014 & 0.014 & 0.015 \\
\hline 
$1 \cdot 10^{-4}$ & 0.05 & 0.034 & 0.033 & 0.035  & 0.038 & 0.04 & 0.044 & 0.044 & 0.047\\
\hline 
$1 \cdot 10^{-3}$ & 0.092 & 0.096 & 0.097 & 0.106 & 0.115 & 0.124 & 0.131 & 0.14 & 0.146 \\
\hline 
\end{tabular}
}
\end{table}

\begin{table}[ht]
\centering
\caption{The condition number $C(\mathbf{A_N,\beta})$ depending on the parameters $N$ and $\beta$ with zero noise, Example 4.1.}
%\footnotesize
\label{table 2}
\scriptsize
\renewcommand{\arraystretch}{1.2} % Increases the height of each row
\begin{tabular}{|c|c|c|c|c|c|c|c|c|c|}
\hline $\beta \backslash N$ & 4 & 6 & 8 & 10 & 12 & 14 & 16 & 18  & 20 \\
\hline
 0 & 266.57 & $2.27 \cdot 10^3$ & $2.45 \cdot 10^4$ & $4.86 \cdot 10^5$ & $1.44 \cdot 10^7$ & $5.25 \cdot 10^8$ & $1.69 \cdot 10^{10}$ & $1.69 \cdot 10^{10}$ & $3.22 \cdot 10^{13}$ \\ 
\hline
$1 \cdot 10^{-13}$  & 266.57 & $2.27 \cdot 10^3$ & $2.45 \cdot 10^4$ & $4.86 \cdot 10^5$ & $1.44 \cdot 10^7$ & $5.25 \cdot 10^8$ & $1.69 \cdot 10^{10}$ & $6.30 \cdot 10^{11}$ & $2.99 \cdot 10^{13}$ \\
\hline
$1 \cdot 10^{-12}$  & 266.57 & $2.27 \cdot 10^3$ & $2.45 \cdot 10^4$ & $4.86 \cdot 10^5$ & $1.44 \cdot 10^7$ & $5.25 \cdot 10^8$ & $1.69 \cdot 10^{10}$ & $6.19 \cdot 10^{11}$ & $1.93 \cdot 10^{13}$ \\
\hline
$1 \cdot 10^{-11}$  & 266.57 & $2.27 \cdot 10^3$ & $2.45 \cdot 10^4$ & $4.86 \cdot 10^5$ & $1.44 \cdot 10^7$ & $5.25 \cdot 10^8$ & $1.70 \cdot 10^{10}$ & $5.26 \cdot 10^{11}$ & $6.28 \cdot 10^{12}$ \\
\hline
$1 \cdot 10^{-10}$  & 266.57 & $2.27 \cdot 10^3$ & $2.45 \cdot 10^4$ & $4.86 \cdot 10^5$ & $1.44 \cdot 10^7$ & $5.26 \cdot 10^8$ & $1.86 \cdot 10^{10}$ & $2.02 \cdot 10^{11}$ & $1.37 \cdot 10^{12}$ \\
\hline
$1 \cdot 10^{-9}$  & 266.57 & $2.27 \cdot 10^3$ & $2.45 \cdot 10^4$ & $4.86 \cdot 10^5$ & $1.44 \cdot 10^7$ & $5.34 \cdot 10^8$ & $1.29 \cdot 10^{11}$ & $5.43 \cdot 10^{11}$ & $7.49 \cdot 10^{10}$ \\
\hline
$1 \cdot 10^{-8}$  & 266.57 & $2.27 \cdot 10^3$ & $2.45 \cdot 10^4$ & $4.87 \cdot 10^5$ & $1.41 \cdot 10^7$ & $6.16 \cdot 10^8$ & $5.48 \cdot 10^9$ & $8.38 \cdot 10^9$ & $3.11 \cdot 10^{10}$ \\
\hline
$1 \cdot 10^{-7}$ & 266.57 & $2.27 \cdot 10^3$ & $2.45 \cdot 10^4$ & $4.93 \cdot 10^5$ & $1.21 \cdot 10^7$ & $5.85 \cdot 10^8$ & $4.48 \cdot 10^8$ & $5.82 \cdot 10^8$ & $3.83 \cdot 10^9$ \\
\hline
$1 \cdot 10^{-6}$ & 266.54 & $2.27 \cdot 10^3$ & $2.42 \cdot 10^4$ & $5.61 \cdot 10^5$ & $4.65 \cdot 10^6$ & $7.61 \cdot 10^7$ & $6.08 \cdot 10^8$ & $8.14 \cdot 10^8$ & $1.18 \cdot 10^{10}$  \\
\hline
$1 \cdot 10^{-5}$  & 266.32 & $2.30 \cdot 10^3$ & $2.17 \cdot 10^4$ & $2.15 \cdot 10^6$ & $2.39 \cdot 10^6$ & $7.03 \cdot 10^6$ & $2.28 \cdot 10^7$ & $3.24 \cdot 10^7$ & $8.87 \cdot 10^7$ \\
\hline
$1 \cdot 10^{-4}$  & 264.09 & $2.60 \cdot 10^3$ & $1.07 \cdot 10^4$ & $1.98 \cdot 10^5$ & $1.81 \cdot 10^5$ & $7.60 \cdot 10^5$ & $3.36 \cdot 10^7$ & $3.80 \cdot 10^6$ & $5.66 \cdot 10^6$ \\
\hline
$1 \cdot 10^{-3}$ & 243.72 & $8.39 \cdot 10^3$ & $2.81 \cdot 10^3$ & $1.70 \cdot 10^4$  & $2.30 \cdot 10^4$ & $7.12 \cdot 10^4$ & $9.33 \cdot 10^4$ & $8.75 \cdot 10^4$ & $1.20 \cdot 10^5$ \\
\hline
\end{tabular}
\end{table}

\section{\textbf{Numerical Results and Discussion}}

In this section, we explore two examples focusing on inverse Stefan problems in one spatial dimension, where the coefficients of the heat equation are set as $a=1$, and those in the boundary conditions are $\lambda=1$, $L=1$, $\gamma=1$ and $T=1$.

To quantify the deviation of the numerical solution from the exact solution, we introduce the root mean square relative error, denoted as \(\Delta u\) for the heat distribution and \(\Delta P\) for the heat flux:
\begin{equation}
    \Delta u = \left[ \frac{\int_{0}^{T} \int_{0}^{s(t)} \left( u(x,t) - u_N(x,t) \right)^2 dx \, dt}{\int_{0}^{T} \int_{0}^{s(t)} \left( u(x,t) \right)^2 dx \, dt} \right]^{1/2}, 
    \quad
    \Delta P = \left[ \frac{\int_{0}^{T} \left( -\lambda \frac{\partial u_N(0, t)}{\partial x} - P(t) \right)^2 dt}{\int_{0}^{T} \left( P(t) \right)^2 dt} \right]^{1/2}.
\end{equation} Here, \(u(x,t)\) and \(P(t)\) represent the exact heat temperature and heat flux, while \(u_N(x,t)\) and \(\frac{\partial u_N(0, t)}{\partial x}\) denote the reconstructed heat temperature and heat flux.

\begin{table}[ht]
\caption{The values $\Delta P$ of the relative error depending on the parameters $N$ and $\beta$ for different noise levels $\epsilon$, Example 4.1.}
\centering
\label{table 3}
\scriptsize
\renewcommand{\arraystretch}{1.2} % Increases the height of each row
\resizebox{\columnwidth}{!}{
\begin{tabular}{|c|c|c|c|c|c|c|c|c|}
\hline &  \multicolumn{4}{ c| }{$N=6$}   &\multicolumn{4}{c|}{$N=8$}\\
\hline $\beta$  & $C(\mathbf{A_N,\beta})$  & $\Delta P_{\epsilon=1\%}$ & $\Delta P_{\epsilon=3\%}$ & $\Delta P_{\epsilon=5\%}$ & $C(\mathbf{A_N,\beta})$ & $\Delta P_{\epsilon=1\%}$ & $\Delta P_{\epsilon=3\%}$ & $\Delta P_{\epsilon=5\%}$  \\
\hline
0 & $2.27 \cdot 10^3$ & $1.19 \cdot 10^{-3}$ & $7.63 \cdot 10^{-4}$ & $1.38 \cdot 10^{-3}$ & $2.45 \cdot 10^4$ & $6.13 \cdot 10^{-4}$ & $4.14 \cdot 10^{-4}$ & $5.67 \cdot 10^{-4}$ \\
\hline 
$1 \cdot 10^{-7}$ & $2.27 \cdot 10^3$ & $2.07 \cdot 10^{-3}$ & $1.45 \cdot 10^{-3}$ & $2.23 \cdot 10^{-3}$ & $2.45 \cdot 10^4$ & $1.44 \cdot 10^{-3}$ & $1.90 \cdot 10^{-3}$ & $1.67 \cdot 10^{-3}$ \\
\hline 
$1 \cdot 10^{-6}$ & $2.27 \cdot 10^3$ & $5.00 \cdot 10^{-3}$ & $4.44 \cdot 10^{-3}$ & $5.05 \cdot 10^{-3}$ & $2.42 \cdot 10^4$ & $4.45 \cdot 10^{-3}$ & $4.89 \cdot 10^{-3}$ & $4.73 \cdot 10^{-3}$ \\
\hline 
$1 \cdot 10^{-5}$ & $2.30 \cdot 10^3$ & 0.01 & 0.01 & 0.01 & $2.17 \cdot 10^4$ & 0.01 & 0.01 & 0.01 \\
\hline 
$1 \cdot 10^{-4}$ & $2.60 \cdot 10^3$ & 0.03 & 0.03 & 0.03 & $1.07 \cdot 10^4$ & 0.03 & 0.03 & 0.03 \\
\hline 
$1 \cdot 10^{-3}$ & $8.39 \cdot 10^3$ & 0.10 & 0.10 & 0.10 & $2.81 \cdot 10^3$ & 0.10 & 0.10 & 0.10 \\
\hline \hline &  \multicolumn{4}{ c| }{$N=10$}   &\multicolumn{4}{c|}{$N=16$}\\
\hline $\beta$  & $C(\mathbf{A_N,\beta})$  & $\Delta P_{\epsilon=1\%}$ & $\Delta P_{\epsilon=3\%}$ & $\Delta P_{\epsilon=5\%}$ & $C(\mathbf{A_N,\beta})$ & $\Delta P_{\epsilon=1\%}$ & $\Delta P_{\epsilon=3\%}$ & $\Delta P_{\epsilon=5\%}$  \\
\hline
0 & $4.86 \cdot 10^5$ & $1.11 \cdot 10^{-3}$ & $9.18 \cdot 10^{-4}$ & $2.24 \cdot 10^{-3}$ & $1.69 \cdot 10^{10}$ & $9.84 \cdot 10^{-3}$ & 0.02 & 0.02 \\
\hline 
$1 \cdot 10^{-7}$ & $4.93 \cdot 10^5$ & $1.68 \cdot 10^{-3}$ & $4.75 \cdot 10^{-4}$ & $1.76 \cdot 10^{-3}$ & $4.48 \cdot 10^8$ & $1.94 \cdot 10^{-3}$ & $2.01 \cdot 10^{-3}$ & $2.57 \cdot 10^{-3}$ \\
\hline 
$1 \cdot 10^{-6}$ & $5.61 \cdot 10^5$ & $4.21 \cdot 10^{-3}$ & $3.33 \cdot 10^{-3}$ & $4.01 \cdot 10^{-3}$ & $6.08 \cdot 10^8$ & $5.13 \cdot 10^{-3}$ & $5.18 \cdot 10^{-3}$ & $5.44 \cdot 10^{-3}$ \\
\hline 
$1 \cdot 10^{-5}$ & $2.15 \cdot 10^6$ & 0.01 & 0.01 & 0.01 & $2.28 \cdot 10^7$ & 0.01 & 0.01 & 0.01 \\
\hline 
$1 \cdot 10^{-4}$ & $1.98 \cdot 10^5$ & 0.04 & 0.04 & 0.04 & $3.36 \cdot 10^7$ & 0.04 & 0.04 & 0.04 \\
\hline 
$1 \cdot 10^{-3}$ & $1.70 \cdot 10^4$ & 0.11 & 0.11 & 0.11 & $9.33 \cdot 10^4$ & 0.13 & 0.13 & 0.13 \\
\hline
\end{tabular}
}
\end{table}

To assess the solution's sensitivity to various noise levels in the measured boundary data, we introduce random noise to the original data. For this purpose, we employ a procedure to generate a set of random real-valued data, which we subsequently add to the exact Neumann data (\ref{96}), how it was done in \cite{Samat} as follows:
\begin{equation}u^{\epsilon}_x(s(t),t) = L \gamma s'(t) + M(0,\sigma^2) \label{07}\end{equation} 
In this equation, $M(0,\sigma^2)$ represents a normal distribution with a mean of zero and a standard deviation of $\sigma = \epsilon \cdot |u_x(s(t),t)| = \epsilon$, where $\epsilon$ denotes the relative (percentage) noise level.

\subsection{Example 1}~\par
 \vskip1mm 
To demonstrate the proposed method, we consider a specific benchmark example for which an analytical solution can be derived (see \cite{Johansson1, Murio, Samat}). The exact solution and the moving boundary are given by:

\begin{equation}
    u(x,t) = -1 + \exp \left( 1 - \frac{1}{\sqrt{2}} + \frac{t}{2} - \frac{x}{\sqrt{2}} \right), \quad \text{in} \quad 0 < x < s(t), \quad 0 < t \leq 1, \label{75}
\end{equation}
\begin{equation}
    s(t) = \sqrt{2} - 1 + \frac{t}{\sqrt{2}}.
\end{equation}

The initial and boundary conditions for the problem are given as follows. 
$$
\begin{gathered}
f(x)=-1+\exp \left( 1-\frac{1}{\sqrt{2}}-\frac{x}{\sqrt{2}}\right), \quad s_0=\sqrt{2}-1, \quad x \in[0, s(0)], \\
P(t)=\frac{-1} {\sqrt{2}}\exp \left( 1-\frac{1}{\sqrt{2}}+\frac{t}{2}\right), \quad u^*=0, \quad  t \in[0,1] .
\end{gathered}
$$

In this section, we reformulate the inverse Stefan problem as discussed in Section 2, where the goal is to find a solution satisfying Eqs. (\ref{02})-(\ref{93}) and recover the Neumann boundary condition (\ref{04}):
\begin{equation}
    \frac{\partial u}{\partial x}(0,t) = \frac{-1}{\sqrt{2}} \exp \left( 1 - \frac{1}{\sqrt{2}} + \frac{t}{2} \right), \quad t \in [0,1]. \label{12}
\end{equation}
We approximate the solution \(u(x,t)\) in the form of Eq. (\ref{74}), where \(v_n(x, t)\) is defined by Eq. (\ref{08}). The coefficients \(c_n\) are determined from the system of equations (\ref{10}), which is obtained by satisfying the initial condition (\ref{03}) and boundary conditions (\ref{05}), (\ref{96}). This system is solved using the matrix inversion method.

The relative errors for the reconstructed heat flux \(P(t)\) for different values of \(\beta\) and \(N\) with noise-free data are shown in Table \ref{table 1}. From the table, it can be observed that for small values of \(N\) and \(\beta\), the regularization parameter \(\beta\) has little effect on the recovery results. However, for larger values of \(N \geq 10\), regardless of how small \(\beta\) is chosen, the unregularized solution provides better results. In other words, for noise-free data, setting \(\beta = 0\) may yield the best results.

\begin{figure}[h!]
\centering
\includegraphics[width=16cm]{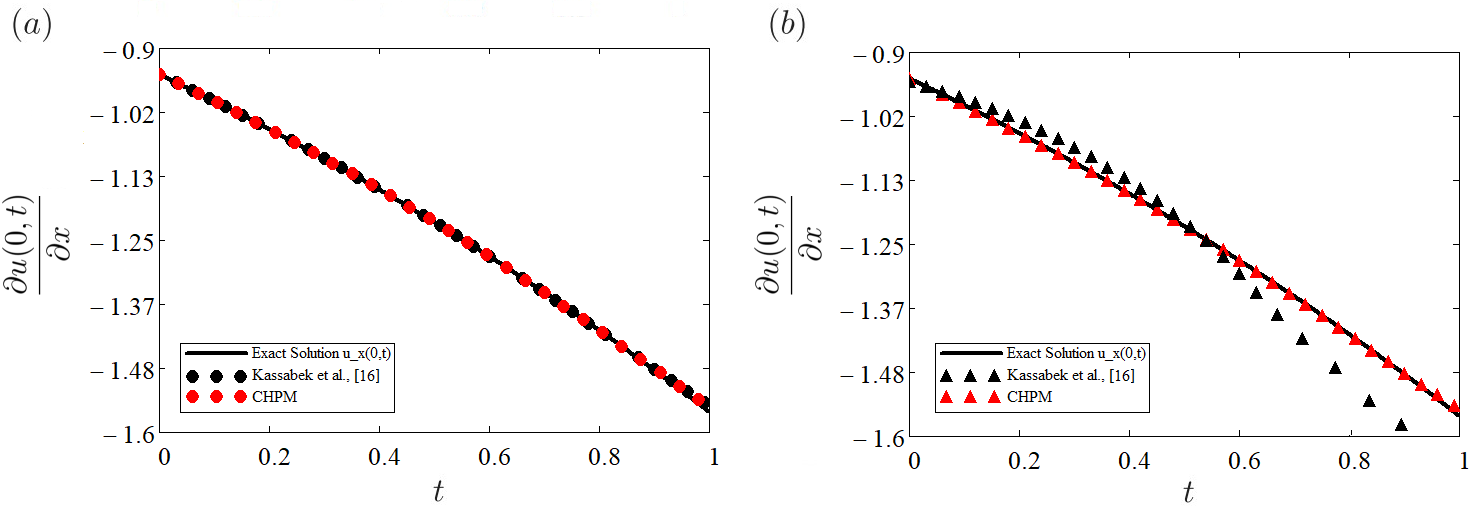}
\caption{For Example 4.1 of the inverse Stefan problem, graphs of the reconstructed heat flux $u_x(0, t)$ for different noise levels: (a) $\epsilon=1\%$, (b) $\epsilon=5\%$ obtained by CHPM with $N=12$ and by Kassabek et al. \cite{Samat} using HPM with $N=12$.}
\label{fig7}
\end{figure}

Since the matrix \(\mathbf{A}\) can be computed independently before any measurements, the condition numbers \(C(\mathbf{A}_N, \alpha)\) can be determined for different values of \(N\) and \(\beta\). This enables the identification of the most suitable combinations of \(N\) and \(\beta\). Table \ref{table 2} shows the condition number \(C(\mathbf{A}_N, \beta)\) for various values of \(\beta\). The results indicate that the condition number is highly sensitive to the choice of \(N\). As seen in the table, \(C(\mathbf{A}_N, \beta)\) increases significantly as \(N\) increases and decreases with larger \(\beta\), revealing a trade-off between stability and accuracy. Stability can be improved by decreasing \(N\) and selecting a sufficiently large \(\beta\), but this comes at the cost of accuracy in approximating the function \(P(t)\).

To assess the sensitivity of the proposed method, we introduce random noise to the boundary data \(s'(t)\) and to \(u_x(s(t),t)\), as outlined in Eq. (\ref{07}). The reconstructed heat flux \(P(t)\) for noise levels of \(\epsilon = 1\%\), \(\epsilon = 3\%\), and \(\epsilon = 5\%\) is shown in Figs. \ref{fig6} and \ref{fig7}, with \(N = 12\). Fig. \ref{fig6}(a) and \ref{fig6}(b) illustrate the CHPM approximations of the Neumann boundary condition (\ref{12}) and the corresponding absolute error for these noise levels, using polynomials of degree \(N = 12\). The numerical results in Fig. \ref{fig7} align well with those presented by Kassabek et al. \cite{Samat}, who used the HPM method with \(N = 12\). It is important to note that the results presented in Figs. \ref{fig6} and \ref{fig7} were obtained without applying regularization, whereas in \cite{Samat}, Tikhonov regularization was used to obtain the numerical results.

Table \ref{table 3} examines various combinations of \((N, \beta)\) with noisy data and identifies the optimal ones. The results show that higher noise levels require a smaller regularization parameter and a reduced cut-off number \(N\). The first test problem considers the time interval \([0, 1]\). To evaluate the method's effectiveness for long-term processes, we also consider different values of \(T\). The numerical results, summarized in Table \ref{table 4}, indicate that increasing \(T\) affects the solution's accuracy.

The numerical results for the first test problem differ from the previous results presented in \cite{Samat}, as shown in Table \ref{table 9}. This table presents the results obtained using two methods: the collocation heat polynomials method (CHPM) and the variational heat polynomials method (VHPM). As seen in Table \ref{table 9}, the VHPM method performs better for \(N \leq 10\), but its accuracy decreases for larger values of \(N\). In contrast, the results obtained by CHPM remain more stable as \(N\) increases. Additionally, the CHPM results are sensitive to the choice of discretization for both the spatial domain \([0, s(0)]\) and the time domain \([0, T]\) when divided into finite intervals.

By integrating conditions (\ref{03}), (\ref{05}), and (\ref{96}) and substituting the approximate solution in the form of (\ref{74}) on the left-hand side, a system of equations is formed. Since there are two boundary conditions on \(s(t)\), namely the Dirichlet condition (\ref{05}) and the Stefan condition (\ref{96}), various discretization strategies can be applied to the time domain \([0, T]\).

\begin{table}[ht]
\caption{Values $P(t)$ of the relative error for different time intervals \([0, T]\) with zero noise, Example 4.1.} \centering
\scriptsize
\renewcommand{\arraystretch}{1.2} % Increases the height of each row
\label{table 4}
\resizebox{\columnwidth}{!}{
\begin{tabular}{|c|c|c|c|c|c|c|c|c|c|}
\hline $T \backslash N$  & 4 & 6 & 8 & 10 & 12 & 14 & 16 & 18  & 20 \\
\hline
1 & $1.84 \cdot 10^{-2}$ & $1.90 \cdot 10^{-3}$ & $1.08 \cdot 10^{-4}$ & $4.80 \cdot 10^{-6}$ & $1.72 \cdot 10^{-7}$ & $5.06 \cdot 10^{-9}$ & $6.57 \cdot 10^{-11}$ & $1.67 \cdot 10^{-11}$ & $2.05 \cdot 10^{-11}$ \\
\hline 
2 & 0.11 & $1.30 \cdot 10^{-2}$ & $1.32 \cdot 10^{-3}$ & $9.48 \cdot 10^{-5}$ & $4.93 \cdot 10^{-6}$ & $4.30 \cdot 10^{-7}$ & $5.06 \cdot 10^{-8}$ & $1.59 \cdot 10^{-9}$ & $1.91 \cdot 10^{-10}$ \\
\hline 
3 & 0.85 & $3.66 \cdot 10^{-2}$ & $4.51 \cdot 10^{-3}$ & $5.23 \cdot 10^{-4}$ & $1.80 \cdot 10^{-4}$ & $4.30 \cdot 10^{-5}$ & $6.32 \cdot 10^{-7}$ & $8.69 \cdot 10^{-7}$ & $3.29 \cdot 10^{-9}$ \\
\hline 
4 & 1.90 & $6.74 \cdot 10^{-2}$ & $1.02 \cdot 10^{-2}$ & $8.17 \cdot 10^{-3}$ & $4.39 \cdot 10^{-3}$ & $2.28 \cdot 10^{-3}$ & $1.62 \cdot 10^{-5}$ & $3.89 \cdot 10^{-5}$ & $2.42 \cdot 10^{-7}$ \\
\hline 
5 & 1.04 & $9.11 \cdot 10^{-2}$ & $7.17 \cdot 10^{-2}$ & 0.13 & $8.77 \cdot 10^{-2}$ & $6.56 \cdot 10^{-3}$ & $6.37 \cdot 10^{-4}$ & $9.94 \cdot 10^{-5}$ & $1.34 \cdot 10^{-5}$ \\
\hline 
\end{tabular}
}
\end{table}

%\tracingtabular
\begin{table}[ht]
\caption{The values of relative error for heat distribution $u(x,t)$ and heat flux $P(t)$ obtained by CHPM and VHPM with zero noise, Example 4.1.} % title of Table
\centering % used for centering table
%\scriptsize
\footnotesize
\renewcommand{\arraystretch}{1.2} 
\begin{tabular}{|c|c|c|c|c|}
\hline 
& \multicolumn{2}{c}{Collocation method} &  \multicolumn{2}{|c|}{Variational method} \\% [1ex] adds vertical space
\hline
$N$  & $ \Delta u $ & $ \Delta P $  & $ \Delta u $ & $\Delta P$      \\ 
\hline 
4 & $0.014$ & $0.018$ & $4.250\cdot 10^{-3}$  & $8.195\cdot 10^{-3}$ \\  % adds vertical space
\hline
6 & $5.064\cdot 10^{-3}$ & $3.064\cdot 10^{-4}$ &  $2.846\cdot 10^{-4}$ & $6.801\cdot 10^{-4}$ \\
\hline
8 & $3.013\cdot 10^{-4}$  & $1.086\cdot 10^{-4}$ & $1.458\cdot 10^{-5}$  & $4.243\cdot 10^{-5}$ \\
\hline
10 & $1.266\cdot 10^{-6}$ & $1.459\cdot 10^{-5}$ & $5.824\cdot 10^{-7}$  & $2.009\cdot 10^{-6}$ \\
\hline
12 & $3.69\cdot 10^{-8}$ & $3.639\cdot 10^{-7}$ & $1.187\cdot 10^{-7}$  &  $7.923\cdot 10^{-8}$ \\
\hline
14 & $9.088\cdot 10^{-10}$ & $5.06\cdot 10^{-9}$ & $2.504\cdot 10^{-7}$  &  $1.190\cdot 10^{-7}$ \\
\hline
\end{tabular}
\label{table 9} % is used to refer this table in the text
\end{table}

An interesting observation regarding the choice of discretization is that the best results are obtained when the spatial interval \([0, s(0)]\) is divided into two subintervals or not divided at all. For example, if we choose to divide \([0, s(0)]\) into two equal subintervals, integrating the initial condition (\ref{03}) over each subinterval, we obtain two equations. To solve a system of size \(N \times N\), the time interval \([0, T]\) for each boundary condition (\ref{05}) and (\ref{96}) can either be divided equally, or for condition (\ref{05}), the number of points used to divide \([0, T]\) must be more than that for the Stefan condition (\ref{96}).

The results shown in Figs. \ref{fig6} and \ref{fig7} were obtained for \(N = 12\), where we solved a \(13 \times 13\) system. In this case, condition (\ref{05}) was divided into 6 intervals, condition (\ref{96}) into 5 intervals, and condition (\ref{03}) into 2 intervals. By integrating these conditions over their respective intervals, we obtain the necessary system of equations.

%\newpage
\subsection{Example 2}~\par
 \vskip1mm
As the second example \cite{Samat}, we consider the moving boundary function
\begin{equation} s(t)=2 a \sqrt{t+t_0}, \end{equation} where $t_0=0.162558$ and $\alpha=0.620063$ with the exact solution given by  \begin{equation} u(x,t)=1-\frac{ \operatorname{erf}\left( \frac{x}{2\sqrt{t+t_0}} \right)}{\operatorname{erf}(\alpha)}.\label{76} \end{equation} 

The initial and boundary conditions considered in the problem are given by
$$
\begin{gathered}
f(x)=1-\frac{ \text{erf}\left( \frac{x}{2\sqrt{t_0}} \right)}{ \text{erf}(\alpha)}, \quad s_0=2 a \sqrt{t_0}, \quad x \in[0, s(0)], \\
P(t)=-\frac{1}{\sqrt{\pi(t+t_0)}\text{erf}(\alpha)}, \quad u^*=0, \quad  t \in[0,1] .
\end{gathered}
$$

In the second test example, the heat flux to be recovered follows the form: \begin{equation}\frac{\partial u}{\partial x}(0,t)=-\frac{1}{\sqrt{\pi(t+t_0)}\text{erf}(\alpha)}, 
\quad t \in \left[0,1\right]. \label{13}\end{equation}

\begin{table}[ht]
\caption{The values $\Delta P$ of the relative error depending on the parameters $N$ and $\beta$ with zero noise, Example 4.2.}\centering
\scriptsize
\label{table 5}
\renewcommand{\arraystretch}{1.2} % Increases the height of each row
\begin{tabular}{|c|c|c|c|c|c|c|c|c|c|c|}
\hline $\beta \backslash N$  & 4 & 6 & 8 & 10 & 12 & 14 & 16 & 18 & 20 \\
\hline 0 & 0.042 & 0.042 & 0.024 & 0.028 & 0.03 & 0.111 & 0.361 & 0.148 & 0.147 \\
\hline $1 \cdot 10^{-13}$  & 0.098 & 0.042 & 0.024 & 0.029 & 0.022 & 0.09 & 0.105 & 0.099 & 0.09 \\
\hline $1 \cdot 10^{-12}$ & 0.098 & 0.042 & 0.024 & 0.032 & 0.013 & 0.125 & 0.174 & 0.17 & 0.189 \\
\hline $1 \cdot 10^{-11}$  & 0.098 & 0.042 & 0.024 & 0.034 & 0.042 & 0.303 & 0.335 & 0.32 & 0.32 \\
\hline $1 \cdot 10^{-10}$ & 0.098 & 0.042 & 0.024 & 0.024 & 0.024 & 0.16 & 0.161 & 0.147 & 0.137 \\
\hline $1 \cdot 10^{-9}$ & 0.098 & 0.042 & 0.05 & 0.053 & 0.09 & 0.123 & 0.141 & 0.162 & 0.176 \\
\hline $1 \cdot 10^{-8}$ & 0.098 & 0.043 & 0.069 & 0.079 & 0.091 & 0.298 & 0.3 & 0.31 & 0.312 \\
\hline $1 \cdot 10^{-7}$ & 0.099 & 0.058 & 0.069 & 0.047 & 0.044 & 0.143 & 0.138 & 0.13 & 0.124 \\
\hline $1 \cdot 10^{-6}$ & 0.099 & 0.131 & 0.139 & 0.139 & 0.135 & 0.151 & 0.153 & 0.174 & 0.185 \\
\hline $1 \cdot 10^{-5}$ & 0.101 & 0.128 & 0.127 & 0.122 & 0.118 & 0.201 & 0.202 & 0.207 & 0.209 \\
\hline $1 \cdot 10^{-4}$  & 0.146 & 0.151 & 0.135 & 0.147 & 0.147 & 0.104 & 0.104 & 0.11 & 0.114 \\
\hline $1 \cdot 10^{-3}$ & 0.264 & 0.29 & 0.239 & 0.239 & 0.236 & 0.225 & 0.226 & 0.226 & 0.226 \\
\hline
\end{tabular}
\end{table}

\begin{table}[ht]
\caption{The condition number $C(\mathbf{A_N,\beta})$ depending on the parameters $N$ and $\beta$  with zero noise, Example 4.2.}
\centering
%\footnotesize
\scriptsize
\label{table 6}
\renewcommand{\arraystretch}{1.2} % Increases the height of each row
\resizebox{\columnwidth}{!}{
\begin{tabular}{|c|c|c|c|c|c|c|c|c|c|}
\hline $\beta \backslash N$ & 4 & 6 & 8 & 10 & 12 & 14 & 16 & 18  & 20 \\
\hline
0 & 134.90 & $2.30 \cdot 10^3$ & $4.26 \cdot 10^4$ & $4.90 \cdot 10^5$ & $1.97 \cdot 10^7$ & $7.57 \cdot 10^8$ & $2.53 \cdot 10^{10}$ & $2.53 \cdot 10^{12}$ & $7.38 \cdot 10^{13}$ \\
\hline
$1 \cdot 10^{-13}$ & 134.90 & $2.30 \cdot 10^3$ & $4.26 \cdot 10^4$ & $4.90 \cdot 10^5$ & $1.97 \cdot 10^7$ & $7.57 \cdot 10^8$ & $2.53 \cdot 10^{10}$ & $1.41 \cdot 10^{12}$ & $7.86 \cdot 10^{13}$ \\
\hline
$1 \cdot 10^{-12}$ & 134.90 & $2.30 \cdot 10^3$ & $4.26 \cdot 10^4$ & $4.90 \cdot 10^5$ & $1.97 \cdot 10^7$ & $7.57 \cdot 10^8$ & $2.53 \cdot 10^{10}$ & $1.38 \cdot 10^{12}$ & $2.18 \cdot 10^{14}$ \\
\hline
$1 \cdot 10^{-11}$ & 134.90 & $2.30 \cdot 10^3$ & $4.26 \cdot 10^4$ & $4.90 \cdot 10^5$ & $1.97 \cdot 10^7$ & $7.57 \cdot 10^8$ & $2.56 \cdot 10^{10}$ & $1.14 \cdot 10^{12}$ & $1.09 \cdot 10^{13}$ \\
\hline
$1 \cdot 10^{-10}$ & 134.90 & $2.30 \cdot 10^3$ & $4.26 \cdot 10^4$ & $4.90 \cdot 10^5$ & $1.97 \cdot 10^7$ & $7.55 \cdot 10^8$ & $2.88 \cdot 10^{10}$ & $3.88 \cdot 10^{11}$ & $2.90 \cdot 10^{12}$ \\
\hline
$1 \cdot 10^{-9}$ & 134.90 & $2.30 \cdot 10^3$ & $4.26 \cdot 10^4$ & $4.90 \cdot 10^5$ & $1.97 \cdot 10^7$ & $7.39 \cdot 10^8$ & $9.67 \cdot 10^{10}$ & $7.98 \cdot 10^{10}$ & $1.92 \cdot 10^{12}$ \\
\hline
$1 \cdot 10^{-8}$ & 134.90 & $2.30 \cdot 10^3$ & $4.26 \cdot 10^4$ & $4.90 \cdot 10^5$ & $1.98 \cdot 10^7$ & $6.06 \cdot 10^8$ & $1.62 \cdot 10^{10}$ & $1.15 \cdot 10^{10}$ & $1.16 \cdot 10^{11}$ \\
\hline
$1 \cdot 10^{-7}$ & 134.90 & $2.30 \cdot 10^3$ & $4.26 \cdot 10^4$ & $4.89 \cdot 10^5$ & $2.17 \cdot 10^7$ & $1.91 \cdot 10^8$ & $6.39 \cdot 10^8$ & $2.98 \cdot 10^9$ & $5.27 \cdot 10^{10}$ \\
\hline
$1 \cdot 10^{-6}$ & 134.90 & $2.30 \cdot 10^3$ & $4.28 \cdot 10^4$ & $4.77 \cdot 10^5$ & $5.50 \cdot 10^7$ & $5.25 \cdot 10^7$ & $1.66 \cdot 10^7$ & $2.42 \cdot 10^9$ & $4.00 \cdot 10^9$  \\
\hline
$1 \cdot 10^{-5}$ & 134.92 & $2.30 \cdot 10^3$ & $4.45 \cdot 10^4$ & $3.76 \cdot 10^5$ & $6.42 \cdot 10^6$ & $6.57 \cdot 10^6$ & $9.04 \cdot 10^6$ & $4.41 \cdot 10^7$ & $7.64 \cdot 10^8$  \\
\hline
$1 \cdot 10^{-4}$ & 135.13 & $2.25 \cdot 10^3$ & $7.28 \cdot 10^4$ & $2.07 \cdot 10^5$ & $5.58 \cdot 10^5$ & $4.73 \cdot 10^5$ & $5.85 \cdot 10^5$ & $3.39 \cdot 10^6$ & $3.38 \cdot 10^6$ \\
\hline
$1 \cdot 10^{-3}$ & 137.35 & $1.87 \cdot 10^3$ & $2.69 \cdot 10^4$ & $9.07 \cdot 10^4$ & $4.13 \cdot 10^4$ & $5.29 \cdot 10^4$ & $1.18 \cdot 10^5$ & $1.17 \cdot 10^5$ & $1.12 \cdot 10^5$ \\
\hline
\end{tabular}
}
\end{table}

Table \ref{table 5} presents the relative errors for the reconstructed heat flux \(P(t)\) with various values of \(\beta\) and \(N\) under noise-free conditions. As in Example 4.1, for smaller values of \(N\) and \(\beta\), the regularization parameter \(\beta\) has minimal impact on the recovery results. However, for larger values of \(N \geq 12\), smaller values of \(\beta\) perform better than the unregularized solution, which contrasts with the results from the previous test problems.

Table \ref{table 6} presents the condition number \(C(\mathbf{A}_N, \beta)\) for different values of \(\beta\). The results demonstrate that, similar to Example 4.1, the condition number is highly sensitive to the selection of \(N\). As shown in Table \ref{table 6}, \(C(\mathbf{A}_N, \beta)\) increases significantly with larger \(N\) and decreases with higher values of \(\beta\).

\begin{figure}[h!]
\centering
\includegraphics[width = 16cm]{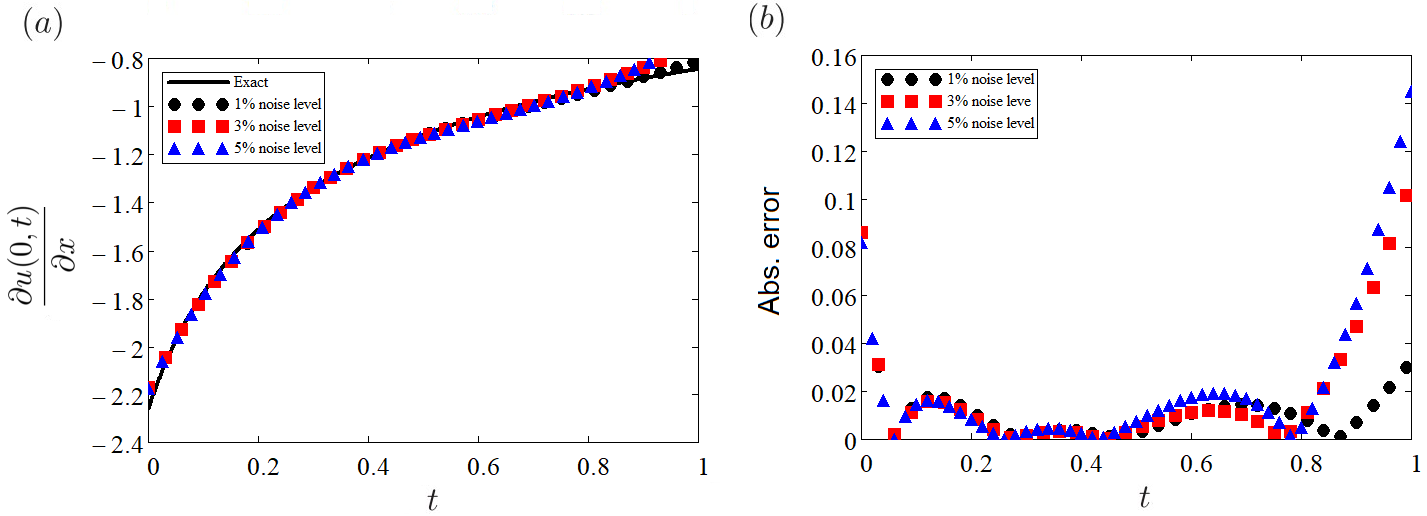}
\caption{For Example 4.2 of the inverse Stefan problem (a) depicts the comparison between the exact and approximate solutions of $u_x(0, t)$, (b) illustrates the absolute error between the exact solution and $u_x(0, t)$ for different noise levels for $N=10$.}
\label{fig10}
\end{figure}

\begin{figure}[h!]
\centering
\includegraphics[width = 16cm]{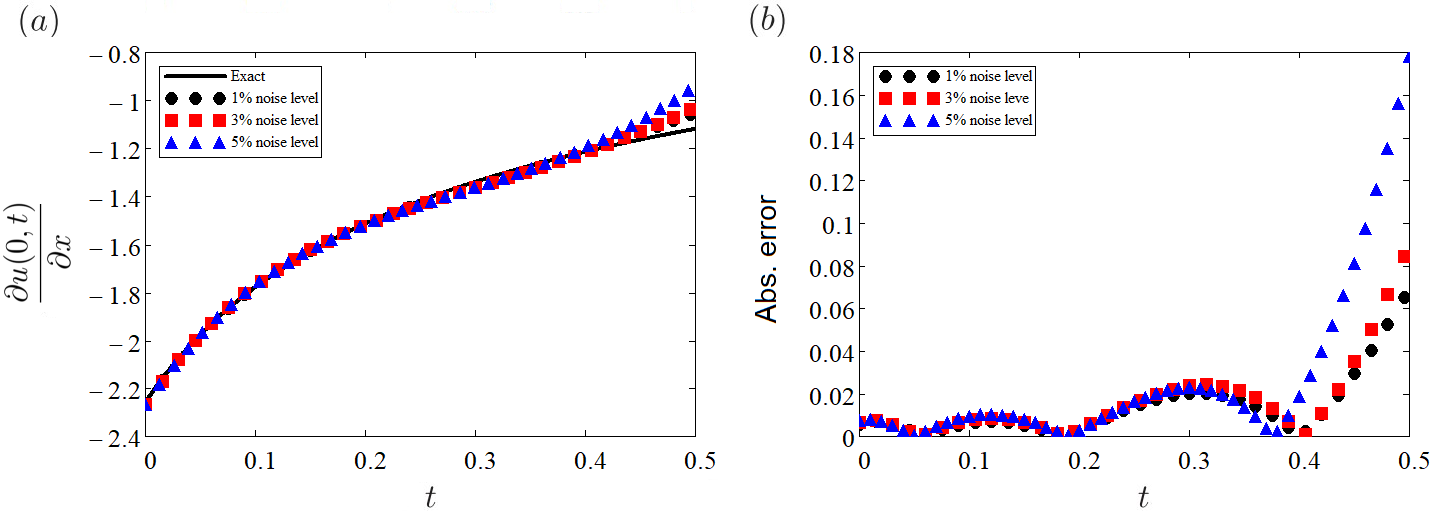}
\caption{For Example 4.2 of the inverse Stefan problem (a) depicts the comparison between the exact and approximate solutions of $u_x(0, t)$, (b) illustrates the absolute error between the exact solution and $u_x(0, t)$ for different noise levels for $N=10$.}
\label{fig11}
\end{figure}

\begin{figure}[h!]
\centering
\includegraphics[width = 16cm]{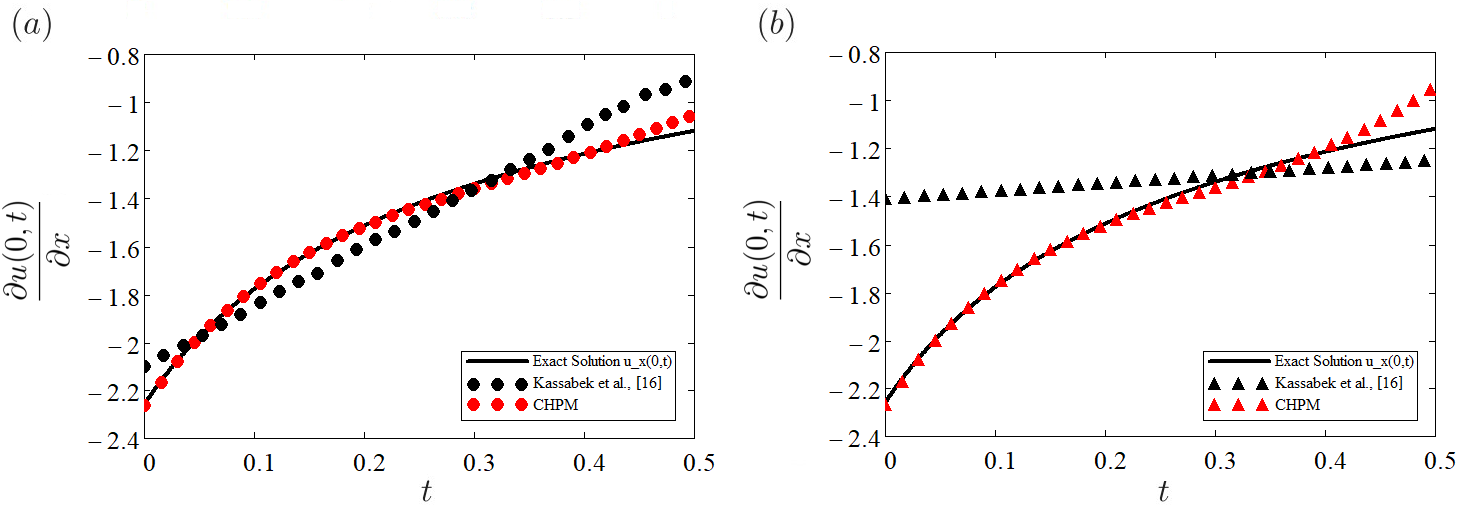}
\caption{For Example 4.2 of the inverse Stefan problem, graphs of the reconstructed heat flux $u_x(0, t)$ for different noise levels: (a) $\epsilon=1\%$, (b) $\epsilon=5\%$ obtained by CHPM with $N=10$ and by Kassabek et al. \cite{Samat} using HPM with $N=10$.}
\label{fig12}
\end{figure}

The reconstructed heat flux \(P(t)\) and the corresponding absolute errors for noise levels of \(\epsilon = 1\%\), \(\epsilon = 3\%\), and \(\epsilon = 5\%\) are shown in Figs. \ref{fig10}, \ref{fig11}, and \ref{fig12} for \(N = 10\). To compare our numerical solution with the previously obtained results by \cite{Samat} (HPM for \(N = 10\)), we consider the time interval \([0, 0.5]\), as illustrated in Fig. \ref{fig12}. It is important to note that no regularization was applied to the results presented in any of the figures.

By examining various combinations of \((N, \beta)\) with noisy data, similar to Example 4.1, we observe that higher noise levels require a smaller regularization parameter and a lower cut-off number \(N\), as shown in Table \ref{table 7}. Additionally, Table \ref{table 8} demonstrates that increasing the time interval \(T\) impacts the accuracy of the solution.

In the case of the second test problem, the CHPM method yielded better results than the VHPM method for all investigated orders of \(N\). The numerical solution results are presented in Table~\ref{table 10}. As in Example~4.1, we employed a non-equal partitioning of the interval \((0, s(0))\) for conditions~(\ref{05}) and~(\ref{96}). The results shown in Figures~\ref{fig10}, \ref{fig11}, and \ref{fig12} were obtained with \(N = 10\), where an \(11 \times 11\) system was solved. Specifically, condition~(\ref{05}) was divided into six intervals, condition~(\ref{96}) into four intervals, and condition~(\ref{03}) into one interval.

\begin{table}[ht]
\caption{The values $\Delta P$ of the relative error depending on the parameters $N$ and $\beta$ for different noise levels $\epsilon$, Example 4.2.}\centering
\scriptsize
\label{table 7}
\renewcommand{\arraystretch}{1.2} % Increases the height of each row
\resizebox{\columnwidth}{!}{
\begin{tabular}{|c|c|c|c|c|c|c|c|c|}
\hline &  \multicolumn{4}{ c| }{$N=6$}   &\multicolumn{4}{c|}{$N=8$}\\
\hline $\beta$  & $C(\mathbf{A_N,\beta})$  & $\Delta P_{\epsilon=1\%}$ & $\Delta P_{\epsilon=3\%}$ & $\Delta P_{\epsilon=5\%}$ & $C(\mathbf{A_N,\beta})$ & $\Delta P_{\epsilon=1\%}$ & $\Delta P_{\epsilon=3\%}$ & $\Delta P_{\epsilon=5\%}$  \\
\hline
0 & $2.308 \cdot 10^{3}$ & 0.041 & 0.0416 & 0.0428 & $4.262 \cdot 10^{4}$ & 0.0228 & 0.0232 & 0.031 \\
$1 \cdot 10^{-7}$ & $2.307 \cdot 10^{3}$ & 0.0561 & 0.0572 & 0.0596 & $4.264 \cdot 10^{4}$ & 0.0672 & 0.0685 & 0.0649 \\
$1 \cdot 10^{-6}$ & $2.307 \cdot 10^{3}$ & 0.1294 & 0.1305 & 0.1333 & $4.28 \cdot 10^{4}$ & 0.1366 & 0.1269 & 0.1366 \\
$1 \cdot 10^{-5}$ & $2.302 \cdot 10^{3}$ & 0.1271 & 0.1276 & 0.1295 & $4.453 \cdot 10^{4}$ & 0.126 & 0.127 & 0.126 \\
$1 \cdot 10^{-4}$ & $2.251 \cdot 10^{3}$ & 0.1517 & 0.1516 & 0.1503 & $7.285 \cdot 10^{4}$ & 0.1352 & 0.1353 & 0.1358 \\
$1 \cdot 10^{-3}$ & $1.871 \cdot 10^{3}$ & 0.2899 & 0.29 & 0.2892 & $2.692 \cdot 10^{4}$ & 0.238 & 0.239 & 0.2388 \\
\hline \hline &  \multicolumn{4}{ c| }{$N=10$}   &\multicolumn{4}{c|}{$N=16$}\\
\hline$\beta$  & $C(\mathbf{A_N,\beta})$  & $\Delta P_{\epsilon=1\%}$ & $\Delta P_{\epsilon=3\%}$ & $\Delta P_{\epsilon=5\%}$ & $C(\mathbf{A_N,\beta})$ & $\Delta P_{\epsilon=1\%}$ & $\Delta P_{\epsilon=3\%}$ & $\Delta P_{\epsilon=5\%}$  \\
\hline 
0 & $4.908 \cdot 10^{5}$ & 0.0174 & 0.0246 & 0.0254 & $2.536 \cdot 10^{10}$ & 0.0543 & 0.1208 & 0.1253 \\
$1 \cdot 10^{-7}$ & $4.894 \cdot 10^{5}$ & 0.0465 & 0.0468 & 0.0455 & $6.396 \cdot 10^{8}$ & 0.139 & 0.1384 & 0.1389 \\
$1 \cdot 10^{-6}$ & $4.771 \cdot 10^{5}$ & 0.1391 & 0.1397 & 0.138 & $1.668 \cdot 10^{7}$ & 0.1523 & 0.1518 & 0.1518 \\
$1 \cdot 10^{-5}$ & $3.767 \cdot 10^{5}$ & 0.1214 & 0.122 & 0.1211 & $9.041 \cdot 10^{6}$ & 0.2019 & 0.2018 & 0.2023 \\
$1 \cdot 10^{-4}$ & $2.075 \cdot 10^{5}$ & 0.1471 & 0.1467 & 0.1469 & $5.856 \cdot 10^{5}$ & 0.1041 & 0.1042 & 0.1038 \\
$1 \cdot 10^{-3}$ & $9.077 \cdot 10^{4}$ & 0.2389 & 0.2387 & 0.2387 & $1.183 \cdot 10^{5}$ & 0.2255 & 0.2257 & 0.2253 \\
\hline
\end{tabular}
}
\end{table}

\begin{table}[ht]
\caption{Values $P(t)$ of the relative error for different time intervals \([0, T]\)  with zero noise, Example 4.2.}\centering
%\scriptsize
\footnotesize
\label{table 8}
\renewcommand{\arraystretch}{1.2} % Increases the height of each row
\resizebox{\columnwidth}{!}{
\begin{tabular}{|c|c|c|c|c|c|c|c|c|c|}
\hline $T \backslash N$  & 4 & 6 & 8 & 10 & 12 & 14 & 16 & 18  & 20 \\
\hline
1 & 0.1 & 0.04 & 0.02 & 0.03 & 0.03 & 0.11 & 0.36 & 0.15 & 0.15 \\
\hline 
2 & 0.43 & 0.54 & 0.55 & 0.69 & 7.45 & 96.25 & 621.04 & $1.21 \cdot 10^{3}$ & $3.58 \cdot 10^{3}$ \\
\hline 
3 & 1.05 & 2.4 & 4.46 & 2.93 & 114.11 & $1.80 \cdot 10^{3}$ & $1.97 \cdot 10^{4}$ & $6.58 \cdot 10^{4}$ & $3.24 \cdot 10^{5}$ \\
\hline 
4 & 1.84 & 5.97 & 16.11 & 28.32 & 700.17 & $1.28 \cdot 10^{4}$ & $2 \cdot 10^{5}$ & $9.46 \cdot 10^{5}$ & $6.53 \cdot 10^{6}$ \\
\hline 
5 & 2.76 & 11.57 & 40.86 & 114.19 & $2.72 \cdot 10^{3}$ & $5.67 \cdot 10^{4}$ & $1.16 \cdot 10^{6}$ & $7.07 \cdot 10^{6}$ & $6.28 \cdot 10^{7}$\\
\hline 
\end{tabular}
}
\end{table}

Overall, the numerical results in this section demonstrate that CHPM approximations offer both accuracy and stability, comparing favorably to exact solutions and yielding better outcomes than those in \cite{Samat}. While increasing the number of heat polynomials initially improves the results, using too many can result in ill-conditioned problems. Nevertheless, with greater computational power, better outcomes may be achieved for higher values of \(N\). From the two numerical test cases, we observe that an error of less than \(5\%\) can be achieved by using \(N = 4\) by solving a \(5 \times 5\) linear system of equations.

To provide a theoretical foundation for our method, we refer to the works of \cite{Samat, Reemtsen, Widder1}, which are summarized in the following theorem:

\newpage

\begin{theorem}
Let \(u(x,t)\) be the unique solution of the problem (\ref{02})-(\ref{93}) and let 
$$ u_{N}(x, t) = \sum_{n=0}^{N} a_{n}^{(N)} v_{n}(x, t) $$ 
be a sequence of polynomials such that:
\begin{itemize}
    \item[(a)] \(\lim _{N \rightarrow \infty} \| S u_{N} - r \| = 0\)
    \item[(b)] \(\left| a_{n}^{(N)} \right| \leq C \left( \frac{e}{2 n t_{0}} \right)^{n / 2}\) for all \(n\) and for \(N = 0,1, \ldots\), where \(C > 0\) is a constant and \(t_{0} > T\).
\end{itemize}
Then:
\begin{enumerate}
    \item \(\lim _{N \rightarrow \infty} a_{n}^{(N)} = b_{n}\)
    \item \(\lim _{N \rightarrow \infty} \left\| \frac{\partial^k}{\partial x^k} u_{N}(x,t) - \frac{\partial^k}{\partial x^k} u(x,t) \right\|_{D} = 0\) for all \(k = 0,1, \ldots\), \\  where \(u(x,t) = \sum_{n=0}^{\infty} b_{n} v_{n}(x,t)\), \(S\) is the trace operator defined by \\ \(S u = \left( u(s(\cdot), \cdot), u_{x}(s(\cdot), \cdot), u(\cdot, 0) \right)\), and the product norm \(r = \left( g(\cdot), \dot{s}(\cdot), f(\cdot) \right)\), with "$\cdot$" representing the independent variables \(t\) or \(x\).
\end{enumerate}
\end{theorem}

\begin{table}[ht]
\caption{The values of relative error for heat distribution $u(x,t)$ and heat flux $P(t)$ obtained by CHPM and VHPM with zero noise, Example 4.2.} % title of Table {ccccc}
\centering % used for centering table
%\scriptsize
\footnotesize
\renewcommand{\arraystretch}{1.2}
\begin{tabular}{|c|c|c|c|c|}
\hline 
& \multicolumn{2}{c}{Collacation method} &  \multicolumn{2}{|c|}{Variational method} \\% [1ex] adds vertical space
\hline 
$N$&$ \Delta u $ & $ \Delta P $  & $ \Delta u $ &$\Delta P$ \\ % inserts table
%heading
\hline % inserts single horizontal line
4& $0.043$ & $0.042$ & $0.201$  & $0.419$ \\ \hline 
6& $0.017$ & $0.042$ &$0.128$ & $0.368$  \\ \hline 
8& $9.599\cdot 10^{-3}$  & $0.024$ & $0.049$  & $0.192$  \\ \hline 
10& $6.862\cdot 10^{-3}$ & $0.028$ & $0.013$  & $0.033$ \\ \hline 
12& $8.281\cdot 10^{-3}$ & $0.03$ & $0.039$  &  $0.170$ \\ \hline 
14& $0.017$ & $0.111$ & $0.043$  &  $0.227$ \\ % [1ex] adds vertical space %inserts single line
\hline
\end{tabular}
\label{table 10}
% is used to refer this table in the text
\end{table}

%\newpage
\section{\textbf{Conclusion}}

Up until now, the heat polynomials method has been applied to solve both one- and two-phase inverse Stefan problems, where the solution in the form (\ref{74}) satisfies the initial and boundary conditions. The unknown coefficients \(c_n \, (n=1,2,...,N)\) are determined from the variational principle, similar to the Trefftz method. In this paper, we applied the collocation method using heat polynomials as basis functions to estimate the temperature distribution and heat flux at the boundary \(x = 0\) for inverse one-phase Stefan problems.

The two numerical examples presented in this paper demonstrate a notably high level of accuracy and stability compared to the exact solution. As shown in Tables \ref{table 1} and \ref{table 5}, an error of less than \(5\%\) can be achieved using \(N = 4\) by solving a \(5 \times 5\) system of equations. Moreover, Table \ref{table 5} shows that in the second test problem, the accuracy of the solution improves significantly by increasing the order \(N\) to 12. Similarly, in the first problem, accuracy improves as \(N\) increases to 18. Based on these results, we conclude that good results can be obtained by choosing values of \(N\) up to 12. At this level, regularization is not needed, as the results indicate that its impact is minimal or does not improve the solution at all.

The numerical results obtained using CHPM differ significantly from those obtained earlier with VHPM, as shown in \cite{Samat}. The VHPM method performs better for \(N \leq 10\), but its accuracy declines for larger values of \(N\), as seen in Table \ref{table 9}. In contrast, the results obtained using CHPM remain more stable as \(N\) increases.

As demonstrated in Theorem 4.1, better results can be achieved by increasing the order of \(N\); however, this increase also raises the condition number of the system, potentially leading to numerical instability. To address this, future research can focus on obtaining more stable solutions in the form of linear combinations of heat polynomials by implementing several strategies.

First, regularization techniques, such as Tikhonov regularization, can help mitigate the ill-conditioning inherent in inverse Stefan problems by penalizing large coefficient values. Additionally, selecting the optimal order \(N\) of the heat polynomials is crucial for balancing accuracy and stability, thereby avoiding overfitting and instability.

Moreover, applying scaling techniques can prevent numerical instabilities caused by large variations in polynomial magnitudes. Preconditioning the system of equations, as outlined in \cite{Liu2013}, can effectively reduce the condition number and further improve stability. Finally, post-conditioning methods, which adjust the coefficients after solving the system, can refine the solution, enhancing accuracy and reducing residual errors \cite{Liu2013post}.

The mathematical simplicity of CHPM makes it well-suited for both nonlinear inverse Stefan problems and direct Stefan problems. Future work will focus on extending the CHPM to address source inverse Stefan-type problems.

\section{\textbf{Acknowledgements}}

The authors sincerely thank the reviewers for their constructive comments and suggestions, which have significantly improved the quality of this paper.

Special thanks go to Prof. Balgaisha Mukanova for her valuable insights and contributions, which greatly helped in refining the numerical aspects of this manuscript.

%\newpage

\end{document}